\documentclass[12pt]{amsart} 
\usepackage[fontsize=12pt]{scrextend}
\usepackage[utf8]{inputenc}
\usepackage{amsmath}
\usepackage{amsfonts}
\usepackage{amssymb}
\usepackage{tikz-cd}
\usepackage{hyperref}
\usepackage[english]{babel}
\usepackage{amsthm}
\usepackage{mathtools}
\usepackage{enumitem}
\usepackage[letterpaper, left=3cm, right=3cm, top = 2cm, bottom = 2cm ]{geometry}
\theoremstyle{definition} 

\newtheorem{de}{Definition}[section]

\newtheorem{propn}[de]{Proposition}
\newtheorem{coro}[de]{Corollary}
\newtheorem{theorem}[de]{Theorem}
\newtheorem{remark}[de]{Remark}
\newtheorem{defn}[de]{Definition}
\newtheorem{thma}[de]{Theorem }

\newtheorem{question}{Question}
\newtheorem{lemma}[de]{Lemma}

\newtheorem*{acknowledgement}{Acknowledgement}

\newcommand{\sL}{{\mathcal L}}
\newcommand{\sO}{{\mathcal O}}
\newcommand{\sP}{{\mathcal P}}
\newcommand{\bP}{{\mathbb P}}

\newcommand{\sV}{{\mathcal V}}

\DeclareMathOperator{\J}{J}   
\DeclareMathOperator*{\Spec}{Spec}    
\DeclareMathOperator{\Sy}{Sym}      
\DeclareMathOperator{\Hom}{Hom} 
\DeclareMathOperator{\Ho}{H} 
\DeclareMathOperator{\Li}{L} 
\makeatletter

\newcommand{\Rmnum}[1]{\expandafter\@slowromancap\romannumeral #1@}
\makeatother
\everymath{\displaystyle} 
\begin{document}

\title{Jacobians, Anti-affine groups and torsion points}

\author{A.J. Parameswaran, Amith Shastri K.}
\address{School of Mathematics, Tata Institute of Fundamental Research,  
1 Homi Bhabha Road, Colaba, Mumbai, India}
\email{param@math.tifr.res.in}
\address{Chennai Mathematical Institute, H1 SIPCOT IT Park, Siruseri, Kelambakkam 603103, India}
\email{amithsk@cmi.ac.in}

\begin{abstract}
We give criteria for the Jacobian of a singular curve \( X\) with at most ordinary \(n\)-point singularities to be anti-affine. In particular, for the case of curves with single ordinary double point we exhibit a relation with torsion divisors. If the geometric genus of the singular curve is  atleast 3 and the normalization is non-hyperelliptic and non-bielliptic, then except for finitely many cases the Jacobian of \( X\) is anti-affine. Furthermore, if the normalization is a general curve of genus atleast 3 then the Jacobian of \( X\) is always anti-affine.
 \end{abstract}
\maketitle

\section{Introduction} 
The Jacobian of a smooth projective curve, the moduli space of degree zero line bundles,  form an important class of algebraic groups. 
They are group varieties, called as abelian varieties, and together with the theta divisor completely determine the curve.
The Jacobian of a singular curve is also an algebraic group, however it is not a complete variety but is an extension of an abelian variety by an affine algebraic group. 
D'Souza gives a canonical compactification of these Jacobians using rank 1 torsion free sheaves cf. \cite{Cyril}.
In \cite{Bhosle} and \cite{BP} the authors give a desingularization of the compactified Jacobian of a curve with simple nodes. 
The normalization of the compactified Jacobian has a particularly nice description in terms of fibre product of projectivization of rank two bundles defined on the Jacobian of the normalization of the nodal curve.
The Jacobian in particular is embedded as an open set in this fibre product.

A group scheme \( G\)	over a field \(k\) is said to be anti-affine if \( \mathcal{O}_G :=\Ho^0(G, \mathcal{O}_G) = k\), i.e., \( G\) admits no regular functions other than constants.
The anti-affine groups  were first studied by C. Sancho de Salas \cite{Sancho} over algebraically closed fields and he calls them variedades cuasi-abelianas.
 Subsequently, in 2009 Brion \cite{Brion2}, and C. Sancho de Salas and F. Sancho de Salas  \cite{Sanchos} independently generalized to other fields. The methods and terminology are different in both these papers.
 Following Brion  \cite{Brion2}, we use the term {\it anti-affine groups} and in \cite{Sanchos} the same object is referred to as quasi-abelian variety.\par

The easiest examples of anti-affine groups are abelian varieties.
A vector extension of an abelian variety is a group \( G\) which is an  extension of an abelian variety by a unipotent group.
Brion  proves that the universal vector extension of an abelian variety is an anti-affine group and a vector extension is anti-affine if and only if it is the push forward from the universal extension by a surjective map (cf. \cite{Brion2}, Proposition 2.3).
In the same paper, Brion gives a characterization of semiabelian varieties, the torus extensions of abelian varieties, to be anti-affine.
However, the examples of anti-affine semiabelian varieties are hard.  Natural examples of semi-abelian varieties are the Jacobian variety of a nodal curve, i.e., the only singularities are ordinary double points.  
The aim of this paper is to show that anti-affine semi-abelian varieties occur naturally as the Jacobian of some  nodal curves.
The criterion for such Jacobian to be anti-affine is the following theorem.

\begin{thma}\label{A1}
Let \( X\) be a nodal curve with \( n\) ordinary double points and let \( y_1,  \cdots y_n\) be the ordinary double points. \( p: C \to X\) be the desingularization and \( p^{-1}(y_j)= \{x_j,  z_j\} \) be the points that map onto \( y_j\). Then the Jacobian \( \J(X) \) is an anti-affine group if and only if the points \( \{ \sL_{z_j} \otimes \sL_{x_j}^{-1}\}_j \) are linearly independent in \( \J(C)\).
\end{thma}
By linear independence of divisors we mean that the divisors satisfy no linear relation other than the trivial one, cf. Definition \ref{linind}.
Thus, for the case of  a nodal curve with a single ordinary double point, the study of anti-affine Jacobian reduces to the study of torsion divisors of the form \(x-z\).
We note that for a smooth curve of genus \( g \geq 3\) the closure of the set of divisors of the form \(x-z\) forms a singular surface \( im(\phi)\), in \( \J(C)\) (cf. Lemma \ref{C-C}).

The study of the surface \( im(\phi)\) leads us to the main theorem of this paper.
\begin{thma}\label{B1}
 Let \( C\) be a smooth projective curve of genus \(g\geq 3\) and \( C\) is not bielliptic hyperelliptic.
  Then the set of torsion points in \( im(\phi)\) is finite. Further, if \( C\) is a general curve then the set of torsion points is empty.
\end{thma}
 Fix a smooth projective curve \( C\).
We will denote  \( \mathcal{X}\) as the set of nodal cuves (with a single node) with \( C\) as the normalization.
Then the Theorem \ref{B1} asserts that there are only finitely many \( X \in \mathcal{X}\) whose Jacobian \( \J(X)\) is not anti-affine. If \( C\) is a general curve, then for every \(X \in \mathcal{X}, \J(X)\) is anti-affine. 

The cases of curves of  geometric genus less than 3 remains in the Theorem \ref{B1}. However, genus 0 curves are rational  curves and  thus the Jacobian of the rational nodal curve is \( \mathbb{G}_m\) and hence affine. For curves of genera 1 or 2, \( \J(C)\) is generated by the divisors of the form \( x-z \). Thus, there are infinitely many nodal curves \( X\) whose Jacobian \( \J(X)\) are not anti-affine.

We do not have an analogue of Theorem \ref{B1} for nodal curves \( X\) with \( n\) nodes. But for a fixed smooth curve \( C\) of positive genus, we construct a nodal curve \( X\) with \(n\) nodes so that \( \J(X)\) is anti-affine, cf. Theorem \ref{choice}.

\begin{thma}\label{E1}
Let \( C\) be a smooth curve. Then we can identify finitely many pairs of points  so that the Jacobian of the corresponding nodal curve is anti-affine.
\end{thma}
 
If the curve \( X\) has cusps as singularities then \( \J(X)\) is  a unipotent extension of \( \J(C)\). However this is not of interest to this paper, cf. Remark \ref{5.4}. We ask the following question: 
\begin{question}\label{ques}
What are the conditions on the singularities of a curve \( X\) so that \( \J(X)\) is semiabelian?
\end{question}
As answer to the above question we prove the following theorem,
 cf. Proposition \ref{semiabelian}.
\begin{thma}\label{F2}
Let \( X\) be a projective curve with \( y\in X\) be the only singular point which we assume to be  an ordinary \(n\)-point singularity. Then \( \J(X)\) is a semiabelian variety and \( \J(X) = \mathbb{P}(\sP_{x_1} \oplus \cdots \oplus \sP_{x_n}) \setminus (s_{x_1} \bigcup \cdots  \bigcup s_{x_n}) \), where \( s_{x_j}\)s' are the divisors corresponding to the coordinate hyperplanes.
\end{thma} 
The converse to the above result is also an interesting question:
\begin{question}\label{ques2}
 If \( \J(X)\) is semiabelian does the curve \( X\) only have ordinary multipoint singularities? 
 \end{question}
 To answer this question we introduce the closely related notion of seminormality, cf. Definition \ref{KolSemi}.
 
 Initially defined for analytic spaces, Traverso  \cite{Trav} gave a purely algebraic definition of seminormality. Salmon \cite{Sal} proves that for a planar curve a singular point \( y\) is seminormal if and only if \(y\) is an ordinary double point.
Bombieri \cite{Bomb} generalizes this result by dropping the condition on curves being planar and proves that the singular point \(y\) is seminormal if and only if  the curve has  normal crossings at \( y\), i.e., smooth branches with linearly independent tangents at \( y\).
Davis \cite{Davis} further extends this result for rings with '\(k\)-rational normalization'.
We prove that the notion of seminormality of a curve \( X\) is equivalent to the Jacobian \( \J(X)\) being semiabelian, cf. Theorem \ref{seminorm}.
 
\begin{thma}\label{F1}
The following are equivalent:
\begin{enumerate}
\item The curve X is seminormal.
\item The curve X has ordinary singularities. 
\item The Jacobian \(\J(X)\) is semiabelian.
\end{enumerate}

\end{thma}
 
In  \cite{Rosen}, Rosenlicht calls semi-abelian varieties as toroidal groups. We will use the term semiabelian instead of toroidal.  Let \( C\) be a smooth projective curve of positive genus and let \( x_{i,1}, \cdots, x_{i,d_i}\) be \( d_i\) distinct points on \( C\) with each local branch of \( C\) immersed in \(X\) so that they give rise to independent vectors in the Zariski tangent space. 
Let \( X\) be the curve obtained by pinching the points \( x_{i,1}, \cdots, x_{i,d_i}\) to a point \( y_i\). 
The variety of  degree zero line bundles on \( X\) is also  called  the  Jacobian and we will denote it as \( \J(X)\).
  Rosenlicht shows that  \( \J(X) \) is semiabelian  (cf. the discussion after Theorem 3 of \cite{Rosen}). 
  We give a condition for the  Jacobian \( \J(X)\) to be anti-affine, cf. Theorem \ref{6.5}.
\begin{thma}\label{q1}
 Let \( X\) be a projective, singular curve with a multiple ordinary \(d_i\)-point singularity at \(\{ y_i\}_{i= 1,\cdots ,n}\). 
Then the Jacobian \( \J(X) \hookrightarrow \bP(\sO\oplus \sV_{1})\times_{\J(C)} \cdots \times_{\J(C)} \bP(\sO\oplus \sV_{r})\). 
Furthermore \( \J(X)=\bP(\sO\oplus \sV_{1})\setminus \{\mathop{\bigcup}\limits_{i=1,\cdots d_1} s_i\}\times_{\J(C)} \cdots \times_{\J(C)} \bP(\sO\oplus \sV_{r})\setminus \{\mathop{\bigcup}\limits_{i=1,\cdots d_r} s_i\}\), where each \(s_i \) is a divisor  corresponding to a coordinate hyperplane.
Further, \( \J(X)\) is anti-affine if and only if the points   \( \{\sL_{ij}\}_{\substack{i=1,\dots r \\\ j=2,\dots d_i}}\) are linearly independent in \(\J(C)\).
\end{thma}
In the above Theorem, \( \sV_i\) is the vector bundle \( (\sP_{x_{i,2}} \otimes \sP_{x_{i,1}}^{-1})
\oplus \cdots \oplus  (\sP_{x_{i,d_i}} \oplus \sP_{x_{i,1}}^{-1}) \) and we denote the line bundle \( (\sP_{x_{i,j}} \otimes \sP_{x_{i,1}}^{-1})\) by \( \sL_{ij}\) 
 where \( \sP\) denotes the usual Poincar\'e bundle.

As an analogue  of Theorem \ref{E1} we prove the following result.
\begin{thma}\label{last}
 Let \( C\) be a smooth curve which is not hyperelliptic or bielliptic and let the genus \( g(C)\geq 3\). We fix the following data \( (n, d_1, \dots, d_n)\). Then there exists a curve \( X\) with exactly \( n\) singular points, each singular point is an ordinary \( d_i\)-point singularity and the Jacobian  \( \J(X)\) is anti-affine.
 \end{thma}
The structure of the paper is as  follows. In Section \ref{JacNode} we first indicate a proof of the well known fact that the Jacobian of  a nodal curve with one ordinary double point is an extension of an abelian variety by a torus. We will also prove Theorem \ref{A1} for nodal curves with one node. 
 In Section \ref{Prelims}, we study the difference morphism \( \phi\) and show that the variety \( im(\phi)\) is a singular surface. We will also recall the Theorems of Abramovich and Harris, and Raynaud. 
 In Section  \ref{TorsionPt} we prove the Theorem \ref{B1}, by showing that elliptic curves do not occur in the surface \( im(\phi)\). We  use  Kempf's result on the torsion points and the non-emptiness of certain linear systems to  prove the second part of the Theorem \ref{B1}.
In Sec. \ref{MoreNode} we  prove Theorem \ref{A1} for a curve \( X\) with multiple nodes.  We also prove Theorem \ref{E1}, where we construct a nodal curve \( X\) with multiple nodes such that \( \J(X)\) is anti-affine.
In Section \ref{Onp} we first prove Theorem \ref{F2}. We then define seminormality and  give a proof of Theorem \ref{F1}. We then return back to the question of when will the Jacobian \( \J(X)\)  be anti-affine
and prove the criterion Theorem \ref{q1} first for the case of a single ordinary multipoint singularity and then for the general case. We will then finish the section with a proof of Theorem \ref{last}.

\begin{acknowledgement}
The first author was supported by  DAE under PIC 12-R\&D-TFR-5.01-0500, the second author was partially supported by a grant from the Infosys Foundation.
\end{acknowledgement}
 
 \section{Jacobian of a curve with one node}\label{JacNode}

We will show that the Jacobian of a nodal curve is a semiabelian variety by describing the Jacobian of the nodal curve as an open set in a projectivization of rank two bundle, cf. \cite{ Bhosle, BP}. This projective bundle  is the normalization of the compactification described in \cite{Cyril}. 
From this description we will obtain a criterion for \( \J(X)\) to be anti-affine.

Let \( X\) be a nodal curve with ordinary double \( y\).
Let \( p : C \to X\) be the normalization and \( p ^{-1}(y) =\{ x,  z\}\). 
By pullback,  a degree zero line bundle on \( X\) gives a degree zero line bundle on \( C\). 
Thus we have a homomorphism of algebraic groups \( \J(X) \xrightarrow{\pi} \J(C)\) (cf. \cite[Chapter 2,  Exercise 6.8(a)]{Hart}).\par

A description of the Jacobian of a nodal curve:
Let \( \sP\) be the Poincar\'e bundle on \( \J(C) \times C\). 
Let \( \sL_{x} \approx \sP|_{\J(C)\times \{x\}}\) and \( \sL_{z} \approx \sP|_{\J(C)\times \{z\}}\)be the line bundles corresponding to the points \( x\) and \( z\) respectively. 
Define \( \bP:= \bP(\sL_{x}\oplus \sL_{z})\)
Then the Jacobian \( \J(X) \hookrightarrow \bP\).
By \cite[Chapter 2,  Exercise 7.8]{Hart},   \( \bP\) has two sections corresponding to the quotient line bundles \( \sL_{x}\) and \( \sL_{z}\).
Then the complement of image of these sections in \( \bP\) is precisely \( \J(X)\).
By the above description,  it is clear that \( \J(X)\) is an extension of the abelian variety \( \J(C)\) by  a one dimensional torus. 

We will give another way to view  the Jacobian \( \J(X)\) as a semiabelian variety.
By \cite[Chapter 2,  Exercise 6.9(a)]{Hart},  we have a short exact sequence:
\begin{center}
\begin{tikzcd}
0\arrow{r} & \mathop{\oplus}\limits_{P\in X} \tilde{\sO_P^*}/\sO_P^* \arrow{r} & \J(X) \arrow{r}{\pi} & \J(C) \arrow{r} & 0, 
\end{tikzcd}
\end{center}
where for each \( P \in X\), \( \sO_P\) is its local ring and  \( \tilde{\sO_P}\) is the integral closure of \( \sO_P\). We will denote the group of units in a ring \( R\) by  \( R^*\).
Locally the ordinary double point is given by \( \mathbb{C}[X, Y]/(XY) \). 
The normalization of the ordinary double point is given by  \( \tilde{\sO}=\mathbb{C}[t] \times \mathbb{C}[t]\).
We define an isomorphism \( \tilde{\sO^*}/ \sO^* \approx \mathbb{C}^*\) by sending each \( (f, g)\in \sO^*\) to \( f(0)/g(0) \in \mathbb{C}^*\).
Thus the corresponding algebraic group is a torus. 
 
We thus have the following well known result.
\begin{propn}\label{semiabel}
Let \( X\) be a nodal curve. Then the Jacobian \( \J(X)\) is a semiabelian variety. The dimension of the maximal torus of \( \J(X)\) is equal to the number of nodes of \( X\).
\end{propn}

\begin{remark}
Frequently, we will refer to a degree 0 line bundle \(\sL\)  either as line bundle on the curve \(C\) or as a point of the Jacobian \( \J(C)\). 
\end{remark}
We will now give a criterion for a nodal curve with one ordinary double point to be anti-affine.
Let \( X\) be a nodal curve obtained by identifying   the points \( x \neq z\) of a smooth curve  \( C\) of geometric genus \( p_g(C)> 0\). Let \( y \in X\) be the nodal point. 
Then  we have a map of Jacobian: \( \J(X)\xrightarrow{\pi} \J(C)\). 
By Proposition \ref{semiabel},  \( \J(X)\) is a semiabelian variety,  whose maximal torus has rank one.
Choose a base point $x_0\in C$ and consider  \( f: C\to \J(C)\), the Abel-Jacobi map,  given by \( x\mapsto \sO_C(x-x_0) \in \J(C)\).

\begin{theorem}\label{node}
Let \( X\) be a projective curve with one nodal point \( y\) and \( C\) be the normalization of \( X\) and the points \( x, z\in C\) map to \( y\) under the normalization. 
Then the Jacobian \( \J(X)\) of the curve \( X\) is an anti-affine group if and only if the point \( \sL_x \otimes \sL_z^{-1}\) is not torsion in \( \J(C)\).
\end{theorem}
\begin{proof}\renewcommand{\qedsymbol}{\( \blacksquare\)}
We will follow Brion's argument \cite[Proposition 2.1]{Brion2} and compute \( \pi_*\mathcal{O}(\J(X))\) using the description of the group \( \J(X)\) described in the beginning of this section.
 Let \( \mathcal{L}_{x},  \mathcal{L}_{z}\) be the line bundles corresponding to the points \( x \textnormal{ and } z\),  respectively.
   We will denote the tensor product \( \mathcal{L}_{z} \otimes \mathcal{L}_{x}^{-1}\) as \( \mathcal{L} \). 
In the beginning of this section we have described the Jacobian  \( \J(X)\) as follows : \( \J(X)\approx \mathbb{P}(\mathcal{O} \oplus \sL)-\{s_{x}\cup s_{z}\} \),  where \( s_x\) and \( s_z\) are the two sections corresponding to the two quotient line bundles \( \sO\) and \( \sL\) respectively.
The regular functions on \( \J(X) \) are given by,  \( \sO( \J(X) )= \Ho^0( \J(C), ~\pi_*(\sO_{\J(X)}) \). We will give the description of \( \pi_*(\sO_{\J(X)}) \). 
Locally,  by considering affine open sets of \( \J(C) \),  we have that \( \pi_*( \sO_{\J(X)})|_{ \Spec A} \approx \mathop{\oplus}\limits_{i \in \mathbb{Z}} \sL^i|_{\Spec A} \). 
Thus,  we have  \( \pi_*(\sO_{\J(X)}) \approx \mathop{\oplus}\limits_{i \in \mathbb{Z}} \sL^i\),  
and \( \Ho^0 ( \J(C), \pi_*(\sO_{\J(X)}) \approx \Ho^0(\J(C), \mathop{\oplus}\limits_{i \in \mathbb{Z}} \sL^i)
\approx \mathop{\oplus}\limits_{i\in \mathbb{Z}}\Ho^0(\J(C),  \sL^i)\).
Note that,  \( \Ho^0(\J(C),  \sL^0)= \mathbb{C}\). 
Thus,  \( \J(X)\) is anti-affine if and only if \( \Ho^0(\J(C),  \sL^i)=0\) for each \( i \neq 0\). Since 
 \( \sL^i\) are degree zero line bundles, the only line bundle with non zero sections is the trivial line bundle. Thus \( \J(X)\) is anti-affine if and only if the point \( \sL = \sL_x \otimes \sL_z^{-1}\) is not torsion in \( \J(C)\).\par
\end{proof}

\section{The Difference Morphism and the Jacobian}\label{Prelims}
By Theorem \ref{node} to study the property of \( \J(X)\) to be anti-affine, we look at the difference morphism \( \phi : C\times C \to \J(C)\) given by \( (x, y) \mapsto \sO_C(x-y) \). Here \( C\) is the normalization of the nodal curve \( X\).
We will denote the image of this morphism by  \( im( \phi)\). 
We will see that the variety \( im(\phi)\) is a singular surface in \( \J(C)\) when the genus \( g\) of \( C\) is larger than two.
We will  recall a result of Abramovich and Harris, which characterizes smooth curves whose second symmetric product contains an elliptic curve, cf. Proposition \ref{DH}. 
We will also recall a result of Raynaud, which characterizes subvarieties of an abelian variety which have dense set of torsion points, cf. Proposition \ref{Ray}.

Let \( C\) be a hyperelliptic curve and  \( \sigma\) be the hyperelliptic involution on \(C\).
Then we have the following lemma that reduces the study of \(im(\phi)\) to \(\Sy^2(C)\). 

\begin{lemma}\label{hypsym}

 For a hyperelliptic curve \( C\),  the map \( \phi\) factors through \(C^{(2)}:= (C\times C)/ \sim\),  where \( (x, y) \sim (\sigma(y),  \sigma(x))\) and \( C^{(2)}\) is isomorphic to \( \Sy^2 C\).
\end{lemma}
\begin{proof}\renewcommand{\qedsymbol}{\( \blacksquare\)}

Consider the morphism \( id \times \sigma_2: C\times C \to C \times C\). 
This morphism descends to a morphism between \( \Sy^2 C \to C^{(2)} \), (the two members of the equivalence class in \( \Sy^2 C\) are mapped to the same class in \( C^{(2)}\) ) which is an isomorphism since \( \sigma \) is an involution.
Further,  since \( \sigma\) is the hyperelliptic involution on \( C\),  the divisor \( \sigma(y)-\sigma(x)\) has the same image as the divisor \( x-y\) in \( \J\). 
Thus,  \( \phi\) factors through \( C^{(2)}\). 
To put it succinctly,  we have the following commutative diagram:

\begin{center}
\begin{tikzcd}
C \times C \arrow{r}{id \times \sigma} \arrow[swap]{d} & C \times C \arrow{r}{\phi} \arrow[swap]{d} & J\\
\Sy^2 C \arrow{r} & C^{(2)} \arrow{ur}.\\
 
\end{tikzcd}
\end{center}

\end{proof}
Note that the above lemma shows that the morphism \( \phi\) has degree two away from the diagonal.

\begin{lemma}\label{blowdown}
Let \( C\) be a smooth projective curve of genus greater than \( 2\) and non-hyperelliptic.
Then the map \( \phi : C\times C \to \J(C)\) is birational onto the image,  in particular the diagonal blows down to a point.
\end{lemma}
\begin{proof}\renewcommand{\qedsymbol}{\( \blacksquare\)}

Suppose that \( (x_i,  z_i) \in C \times C\setminus \Delta,  i=1, 2\) be two points. 
Further, assume that \( \sL_{z_1}\otimes \sL_{x_1}^{-1}=\sL_{z_2} \otimes \sL_{x_2}^{-1}\),  i.e.,  the divisors \( z_1-x_1\) and \( z_2-x_2\) are linearly equivalent and hence the divisors  \( z_1+x_2,  z_2+x_1\) are also linearly equivalent.
Since \( x_i \neq z_i\) and if the divisors \(  z_1+x_2,  z_2+x_1\) are distinct then we get a map from \( C \to \mathbb{P}^1\) with zeroes at \( z_2,  x_1\) and poles at \( z_1,  x_2\). 
This gives a degree 2  morphism  from \( C \to \mathbb{P}^1\)   and hence \( C\) is a hyperelliptic curve which is contrary to our assumption.
Similarly,  assume \( x_1= z_2 \) and \( x_2 =z_1\),  this also gives a morphism of degree 2 from \( C\to \mathbb{P}^1 \) with a zero of order 2 at \( z_2\) and a pole of order 2 at \( z_1\) which is also a contradiction.\\
Therefore assume that \( x_1=x_2\),  this forces \(z_1=z_2\) and conversely. 
Thus,  \( \phi |_{C \times C \setminus \Delta}\) is an injection. 
\end{proof}

\begin{lemma}\label{C-C}
Let \( C\) be a curve of genus \( g\geq 3\). Then the surface \( im( \phi)\) is singular.
\end{lemma}
\begin{proof}\renewcommand{\qedsymbol}{\( \blacksquare\)}
Let \( C\) be a non-hyperelliptic curve.  Under the map \( \phi\) the diagonal \( \Delta\) blows down to a point and  \( \Delta ^2 =2-2g < -1\). Hence the surface \( im( \phi)\) is singular.\par
In case  \( C\) is hyperelliptic, by Lemma \ref{hypsym} the surface \(im(\phi) \approx \Sy^2(C) \subset \J^2(C)\). 
For the surface \( \Sy^2(C)\) the diagonal \( \Delta\) maps onto \( \mathbb{P}^1\) with self intersection \( 1-g\) (use the adjunction formula for the surface  \( \Sy^2C\) cf.  \cite{Hart} Chapter V Proposition 1.5 ) . 
Thus \( \Delta\) has self intersection \( -1\) only when genus of \( C\) is 2. 
Thus the surface \( im( \phi)\) is singular for \( C\) hyperelliptic of genus \( g \geq 3\).
\end{proof}
For further questions and properties of the surface \( im(\phi)\), we refer to \cite{Welt}.

Let \( C\) be a smooth projective curve of genus \( g \) and \( C_i,  i=1,  2\) be smooth projective curves of genera \( k_i\). 
Let \( \pi_i: C \to C_i\),  \( i=1, 2\) of degree \( h_i \).
We then have the following inequality which gives an upper bound for the genus \( g\) in terms of the genera \( k_i \) and degrees \( h_i\).
\begin{lemma}[\textit{Castelnuovo-Severi Inequality}]\label{CSI}
As above if \( \pi:= ( \pi_1,  \pi_2 ):C \to C_1 \times C_2 \) is birational onto the image then:
\[ g \leq (h_1-1)(h_2-1)+k_1h_1+k_2h_2\]
\end{lemma}
\begin{proof}\renewcommand{\qedsymbol}{\( \blacksquare\)}
We refer to \cite[Theorem 3.5]{Accola} for a proof.
\end{proof}
\begin{defn}
A smooth projective curve \( C\) is said to be \textit{bielliptic} if there is a degree two map from \( C\) to an elliptic curve.
\end{defn}
Similar to a hyperelliptic curve,  a bielliptic curve has an involution \( \tau\).
This involution in general is not unique. 
However by the Castelnuovo-Severi inequality, Lemma \ref{CSI}, we have the following corollary.
\begin{coro}\label{unibi}
Let \( C\) be a bielliptic curve of genus \( g> 5\). Then  bielliptic involution \( \tau \) is unique.
\end{coro}
\begin{proof}\renewcommand{\qedsymbol}{\( \blacksquare\)}
We will assume that there are two distinct maps from \( C\) onto two elliptic curves \( E_1 \textnormal{ and } E_2\) and show that the genus \( g\) of \( C\) is less than 6.
Let \( p_1: C \to E_1 \textnormal{ and } p_2: C \to E_2\) be degree two morphisms from \( C\). 
Applying the Castelnuovo-Severi inequality to the morphism \( (p_1,p_2): C \to E_1 \times E_2\) we get the bound on the genus.
\end{proof}

We will now show that a smooth projective curve that is both bielliptic and hyperelliptic does not have high genus. Specifically:

\begin{coro}\label{bihyp}
A smooth curve \( C \) that is both hyperelliptic and bielliptic has genus at most 3.
\end{coro}
\begin{proof}\renewcommand{\qedsymbol}{\( \blacksquare\)}
Since \( C\) is both hyperelliptic and bielliptic,  we have  maps  \( \pi _1: C\to \mathbb{P}^1\) and \( \pi_2: C\to E\),  where \( E\) is a genus 1 curve and the maps \( \pi_1\) and \( \pi_2 \) have degree 2. 

Now consider \( \pi:= (\pi_1,  \pi_2): C \to \mathbb{P}^1 \times E\). 
Let \( D\) be the normalization of the image of \( C\). 
We claim that \( D, C\) are birational. On the contrary, suppose that the curves \( C ,D\) are not birational. Then the map from \( C \to D\) has a degree \( d> 1\). The first projection gives a map from \( C \to \mathbb{P}^1\) of degree 2 factoring by \( D\). Thus \( D \approx \mathbb{P}^1\). 
Similarly taking the second projection we get that \( D \approx E\). But this is not possible. Thus the map from \( C\to D\) has no degree and  hence \( C\) and \( D\) are birational to each other. 
Now by the Lemma \ref{CSI} we have \( g \leq 3\).
\end{proof}

We will now state a theorem of Abramovich-Harris which characterizes the curves whose second symmetric product contains an elliptic curve.
\begin{propn}\label{DH}
Let \( C\) be a smooth projective curve. If the symmetric product \( \Sy^2(C)\) contains a genus 1 curve, then \( C\) is a bielliptic curve or genus of \( C\) is at most 2. If  \(g(C)> 3\) then \( C\) is not hyperelliptic. 
\end{propn}
\begin{proof}\renewcommand{\qedsymbol}{\( \blacksquare\)}
See \cite[Theorem 3]{DH} for a proof.
\end{proof}
We note that a weaker version of the above result was proved by Harris and Silverman, cf. \cite{HS} Theorem 2.

Moving from curves, we will now state a crucial result of Raynaud which characterizes the subvarieties of an abelian variety which have dense set of torsion points.
\begin{propn}\label{Ray}
Let \( A\) be an abelian variety and \( X\) a subvariety. Let \( A_{Tor}\) be the torsion points of \( A\). If \( A_{Tor} \cap X\) is a dense subset in \( X\),  then \( X\) is a translate of an abelian subvariety \( B \subset A\) by a torsion point of A.
\end{propn}
\begin{proof}\renewcommand{\qedsymbol}{\( \blacksquare\)}
See \cite{Raynaud} for a proof. 
\end{proof}
The case of curves arising as subvarieties of abelian varieties has also been dealt by Raynaud \cite{Ray2}.
\section{Torsion Divisors of the form \(x-z\)}\label{TorsionPt}
We will prove the Theorem \ref{B1} in this section. 
We are interested in the torsion points of the form \( x-z\) i.e., the torsion line bundles in \( \J(C)\) corresponding to the divisor \( x-z\) in the surface \( im(\phi)\). 
By Raynaud's result, to conclude finiteness of torsion points, it is enough to show that \(im(\phi)\) has no elliptic curves.
We discuss this in two subsections: First we deal with the case when \( C\) is non-hyperelliptic. 
In this case the proof is straightforward.
The case when \( C\)  is hyperelliptic is more involved and is dealt subsequently.

We will handle the second part of the Theorem later in this section.

\subsection{\( C\) is not hyperelliptic}

We will first show that for a smooth,  projective curve,  \( C \) of genus \( g \geq 3\) and non-hyperelliptic the image of the morphism \( \phi \) has only finitely many torsion points.

\begin{lemma}\label{elliptic}
Let \( C\) be a smooth projective curve of genus greater than \( 2\) and non-hyperelliptic.
Then the image of the map \( \phi: C\times C \to \J(C)\) has no elliptic curves.
\end{lemma}
\begin{proof}\renewcommand{\qedsymbol}{\( \blacksquare\)}
Suppose that \( i:E \to \phi (C\times C )\) be a curve of genus \( 1\) in the image of \( \phi\). 
By Lemma \ref{blowdown}, the map \( \phi\) is injective outside the diagonal. Thus we can take the reduced inverse image of the elliptic curve which we denote by \( \tilde{E}\).
 The morphism \( \tilde{i}: \tilde{E} \dashrightarrow C\times C\) is rational (see \cite[Chapter 2,  Section 7,  Corollary 7.15]{Hart}) which extends to the entirety of \( \tilde{E} \) (see \cite[Chapter 1,  Section 6,  Proposition 6.8]{Hart}).
Now,  consider the composite of the projection maps \( \pi_1,  \pi_2 : C\times C \to C\) with \( \tilde{i}\). 
This gives a map from \( \tilde{E }\) a curve of genus \( 1\) to \( C\) and by the Riemann-Hurwitz formula the composite map is constant. 
Thus \( E\) is a point. 
\end{proof}

We will now show that the surface \(im(\phi) \subset \J(C)\) has finitely many torsion points.
\begin{lemma}\label{tors}
Let \( C\) be a smooth projective curve of genus greater than \( 2\) and non-hyperelliptic
then \( \phi: C\times C \to \J(C)\) has only finitely many torsion points in \( im(\phi) \).
\end{lemma}
\begin{proof}\renewcommand{\qedsymbol}{\( \blacksquare\)}

Suppose that there are infinitely many torsion points in the image variety \( im(\phi)\).
 Let \(Y\) be  the Zariski closure of the torsion points. 
 By Raynaud's theorem, Proposition \ref{Ray},  \( Y\) is a translate of an abelian variety. 
 Since the variety \( im(\phi)\) is singular, \( im(\phi)\) cannot be a translate of an abelian surface.
 Thus the subvariety \( Y\) is not  the entirety of \( im(\phi)\). 
 Then  \( Y\) should be a translate of an elliptic curve.
This implies that \( Y\) is a curve of genus \( 1\) in \( im(\phi)\). 
 However by Lemma \ref{elliptic} this is not possible.
Therefore \( im(\phi)\) can contain finitely many torsion points. 
\end{proof}

In case \( C\) is a curve of genus \(g=1 \textnormal{ or } 2 \), \( \phi\) is a surjection onto the \( \J(C)\),  and \(\phi ( C \times C) \)  is an abelian variety and thus has infinitely many torsion points.\par

\subsection{ \( C\) is hyperelliptic}
We need to show that the surface \(im(\phi)\) has no elliptic curves. But by Lemma \ref{hypsym} it is enough to show that the variety \( \Sy^2(C)\) has no elliptic curves.

\begin{propn}\label{ell}
If \( C\) is a hyperelliptic curve of genus \( g \geq 4\),  then surface \( im(\phi) \) has no genus 1 curves.
\end{propn}

\begin{proof}\renewcommand{\qedsymbol}{\( \blacksquare\)}
Let \( E \subset im( \phi)\) be a genus 1 curve. 
Note that \( \Sy^2(C)\) is birational onto its image in \( \J^2(C)\),.
 Further, by Lemma \ref{hypsym} the difference map factors through a variety isomorphic to \( \Sy^2(C)\). 
By Proposition \ref{DH}, \( C\) is bielliptic.
Thus by Corollary \ref{bihyp} genus of \( C\) is at most 3.
However, by assumption the genus is atleast 4.
Thus the surface \( im( \phi) \) has no genus 1 curves.
\end{proof}

We are now ready to show that the surface \(im(\phi)\) admits only finitely many torsion points, even in the case of hyperelliptic curves.
\begin{coro}\label{hyp}
For a hyperelliptic curve of genus \( g \geq 4\),  The surface \( im( \phi) \) has finitely many torsion points.
\end{coro}
\begin{proof}\renewcommand{\qedsymbol}{\( \blacksquare\)}
Suppose that the surface \( im( \phi) \) has infinitely many torsion points. Let \( Y\) be the Zariski closure of these torsion points. 
By Proposition \ref{Ray} \( Y\) is a translate of an abelian variety.
The surface \( im( \phi) \) is singular by Lemma \ref{C-C} and hence \( Y\neq im( \phi)\). 
Similarly,  \( Y\) cannot be a genus 1 curve by  Proposition \ref{ell}.
Thus \(Y \) is a point. 

\end{proof}

Now only the case of hyperelliptic curves of genus 3 remains.
\begin{coro}\label{notbiell}
Let \(C\) be a hyperelliptic curve of genus 3 which is not bielliptic. Then \( im(\phi)\) has finitely many torsion points.
\end{coro}
\begin{proof}\renewcommand{\qedsymbol}{\( \blacksquare\)}
Suppose that \( im(\phi)\) has infinitely many torsion points.
Let \( Y\) be the Zariski closure of the set of torsion pints in \( im(\phi)\).
Since \(im(\phi)\) is singular (Lemma \ref{C-C})  and by Raynaud's theorem \( Y\) is an elliptic curve. 
The isomorphism of  Lemma \ref{hypsym} gives an elliptic curve in \( \Sy^2C\).
Thus \( C\) is also a bielliptic curve of genus 3. However, this is in contradiction to our assumption.
Thus \( im(\phi)\) has finitely many torsion points. 
\end{proof}

We will now state a result of Kempf which gives a necessary cohomological condition for a multiple of a divisor to be rationally equivalent to zero.

\begin{propn}\label{Kempf}
 Let C  be  a  curve   with   general  moduli    in characteristic   zero. For  a  divisor  \(D\) such that \( d . D\) is rationally equivalent to zero for some \( d> 0\),   then the cohomology  
 \(\Ho^1(C,(D_0  + D_{\infty})\)  must   be  zero  (equivalently  \(|K - D_0 - D_{\infty}| \)  is    empty).

\end{propn}
\begin{proof}\renewcommand{\qedsymbol}{\( \blacksquare\)}
We refer to the main Theorem in \cite{Kempf} for a proof.
\end{proof}

As a corollary we have the following corollary.
\begin{coro}\label{gencrv}
Let \( C\) be a general curve of genus atleast 3. Then \( im (\phi)\) has no torsion divisors of the form \( x-z\).
\end{coro}
\begin{proof} \renewcommand{\qedsymbol}{\( \blacksquare \)} 
For any pair of distinct points \( x,z \in C \) and \( g \geq 4 \), the degree of  $ K - x -z \geq g $. And hence the linear system \(|K-x-z|\) is non empty. Thus by Kempf's theorem the corresponding divisor \( x-z\) cannot be a torsion point.

Now, let \( C\) be a general curve of genus 3. Then for points \( x,z \in C\) consider the divisor \( D:=x-z\). 
Suppose that \(d. D\) is rationally equivalent to 0 for some \(d>0\). 
Note that the Euler characteristic \( \mathcal{X}(C, (x+z) ) =0\). 
Thus,  \(\Ho^0(C,\mathcal{O}(x+z))=\Ho^1(C,\mathcal{O}(x+z)) \). 
Thus either both  are non zero or both vanish.
However \( \mathcal{O}(x+z)\) is  an effective divisor and thus \( \Ho^0(C, \mathcal{O}(x+z)) \neq 0\).
Hence \( \Ho^1(C, \mathcal{O}(x+z) \neq 0\), a contradiction to Kempf's theorem. 
Thus \( d.D\) is not rationally equivalent to 0 for any \( d> 0\). 
\end{proof}

We will now summarize our results:
\begin{theorem}\label{finale}
Let \( C\) be a smooth projective curve of genus \( g\geq 3\) and \( C\) is not bielliptic hyperelliptic. Then the set of torsion points in \( im(\phi)\) is finite. Further, if \( C\) is a general curve then the set of torsion points in \( im(\phi)\) is empty.
\end{theorem}

\begin{remark}
Curves of genus 3:\par
Let \( C\) be a smooth curve of genus 3 and not hyperelliptic. 
In this case the canonical divisor \(K_C\) gives an embedding of \(C\) into \(\mathbb{P}^2\) as a smooth quartic. 
Project \(\mathbb{P}^2\) onto \(\mathbb{P}^1\) from a point \( P\in C\). 
Restricting this projection to the curve \( C\) gives a map \( C \to \mathbb{P}^1\).
 Thus \(C\) is a trigonal curve.
 Let \(f: C \to \mathbb{P}^1\) be the degree 3 map. We will compute torsion points of order 3 for these curves.
  A totally ramified point \( Q \) of the map \( f \) gives a torsion point of order three.
   The point \( Q\) is also an inflection point on the curve \( C\).
   To get a divisor of the form \(Q-S \) we need two such divisors.
However there are finitely many inflection points. Thus for a general choice of the projecting  point \( P\) the curve \( C\) will have no 3 torsion points of the form \( Q-S\). 

It should be noted that there are special curves which do admit totally ramified points, and hence have torsion points of order 3. A curve which does not have 3-torsion may have a higher order torsion points.
\end{remark}

\begin{remark}
If the curve \( C\) is of genus 1 or 2, then \( im(\phi)\) is an abelian variety and hence admits infinitely many torsion points. However,  if the divisor \( \pi^{-1}(y)\) corresponding to the ordinary double point \( y \in X\)  is not in this set then the Jacobian \( \J(X)\) is anti-affine.
If \( C\) is a bielliptic hyperelliptic curve of genus 3, then \( im(\phi)\) contains an elliptic curve. Thus in this case there are infinitely many torsion points of the form \( x-z\),
\end{remark}

\section{Nodal Curves with more than one ordinary double point}\label{MoreNode}

 The description of the Jacobian of a curve with one ordinary double point in the second section generalizes to a curve with \( n \) ordinary double points as singularities.
It is an extension of an abelian variety by a torus of rank \( n\).

Let \( X\) be a nodal curve with ordinary double points \( \{y_1,  y_2,  \cdots,  y_n\}\).
Let \( p : C \to X\) be the normalization and \( p ^{-1}(y_j) =\{ x_j,  z_j\}\). 
By pullback,  a degree zero line bundle on \( X\) gives a degree zero line bundle on \( C\). 
Thus we have a homomorphism of algebraic groups \( \J(X) \xrightarrow{\pi} \J(C)\) (cf. \cite[Chapter 2,  Exercise 6.8(a)]{Hart}).\par

A description of the Jacobian of a nodal curve:
Let \( \sP\) be the Poincar\'e bundle on \( \J(C) \times C\). 
Let \( \sL_{x_j} \approx \sP|_{\J(C)\times \{x_j\}}\) and \( \sL_{z_j} \approx \sP|_{\J(C)\times \{z_j\}}\)be the line bundles corresponding to the points \( x_j\) and \( z_j\) respectively. 
Define \( \sP_j:= \bP(\sL_{x_j}\oplus \sL_{z_j})\)
Then the Jacobian \( \J(X) \hookrightarrow \sP_1 \times_{\J(C)} \times \sP_2 \times_{\J(C)} \cdots \times_{\J(C)} \sP_n\).
By \cite[Chapter 2,  Exercise 7.8]{Hart},   each \( \sP_j\) has two sections corresponding to the quotient line bundles \( \sL_{x_j}\) and \( \sL_{z_j}\).
Then the complement of these sections in \( \sP_1 \times_{\J(C)} \times \sP_2 \times_{\J(C)} \cdots \times_{\J(C)} \sP_n\) is precisely \( \J(X)\) cf. \cite[Proposition 2.1]{BP} and the subsequent discussion.
By the above description,  its is clear that \( \J(X)\) is an extension of the abelian variety \( \J(C)\) by \( n -\)dimensional torus.

\begin{defn} \label{linind}
Let \( \{\sL_i\}_{i = 1, \cdots ,n}\) be a set of points in a Jacobian of a smooth curve.
 We say that the divisors  \(\{ \sL_i\}\) are linearly independent if \( \mathop{\sum}\limits_{i=1}^n \sL_i^{k_i} =0\) and \( k_i \in \mathbb{Z}\) for each \( i\) then each \( k_i =0\).
\end{defn}
We will now give a criterion for the Jacobian of a nodal curve with multiple nodes to be anti-affine.

\begin{theorem}\label{highnode}
Let \( X\) be a nodal curve with \( n\) nodes and let \( y_1,  \cdots y_n\) be the nodes. \( p: C \to X\) be the desingularization and \( p^{-1}(y_j)= \{x_j,  z_j\} \) be the points that map onto \( y_j\). Then the Jacobian \( \J(X) \) is an anti-affine group if and only if the points \(\{ \sL_{z_j} \otimes \sL_{x_j}^{-1}\}_{j=1}^n \) are linearly independent in \( \J(C)\).
\end{theorem}
\begin{proof} \renewcommand{\qedsymbol}{\( \blacksquare \)}
We have described the Jacobian, \( \J(X) \) of \( X\) as an open subset of  \( P_1 \times_{\J(C)} P_2 \times_{J(C)} \cdots \times_{J(C)} P_n\) given by removing the image of the pullback of the sections \(s_{x_j},  s_{z_j}\) of \( \sP_j \to \J(C) \). 
As before we will describe \( \pi_* \sO_{ \J(X)} \) in terms of the line bundles \( \sL_j:= \sL_{z_j} \otimes \sL_{x_j}^{-1} \).
We have that,   \( \pi_*(\sO_{\J(X)}) \approx \mathop{\oplus}\limits_{\substack{k_i \in \mathbb{Z}\\\ i \in \{1, 2,  \cdots n\}}} \sL^{k_1}_1 \otimes \sL^{k_2}_2 \otimes \cdots \otimes  \sL^{k_n}_n\) and \( \Ho^0(\J(C), \pi_*(\sO_{ \J(X)} )) \approx \Ho^0 \big(\J(C), \mathop{\oplus}\limits_{\substack{k_i \in \mathbb{Z}\\\ i \in \{1, 2,  \cdots n\}}} \sL^{k_1}_1 \otimes \sL^{k_2}_2 \otimes \cdots \otimes  \sL^{k_n}_n\big)\approx \mathop{\oplus}\limits_{\substack{k_i \in \mathbb{Z} \\\ i\in \{1, 2,  \cdots n\}}}\Ho^0(\J(C),  \sL^{k_1}_1 \otimes \sL^{k_2}_2 \otimes \cdots \otimes  \sL^{k_n}_n  )\).
Therefore the group \( \J(X)\) is anti-affine if and only if the points \( \sL^{k_1}_1 \otimes \sL^{k_2}_2 \otimes \cdots \otimes  \sL^{k_n}_n\) are non trivial, 
i.e.,  the set of points \( \{\sL_1, \sL_2 ,\cdots ,\sL_{n}\}\) are linearly independent in \( \J(C)\).
\end{proof}

By the above theorem it is clear that the Jacobian of a nodal curve being anti-affine depends on the divisors \( \{x_i - z_i\}_{i=1, \cdots, n }\) not satisfying any linear relation. Notice that for a rational nodal curve its Jacobian is a torus and hence is always affine. For other curves we have the following:

\begin{theorem}\label{choice}
Let \( C\) be a smooth non-rational curve. Then we can identify finitely many pairs of points  so that the corresponding Jacobian of the nodal curve is anti-affine.
\end{theorem}

\begin{proof}\renewcommand{\qedsymbol}{\( \blacksquare\)}
We prove this by induction on the number of singularities \( n\). 
For \( n=1\), we choose a point \( \sL \in im(\phi)\) which is not torsion. This is possible as the torsion points in
\( \J(C)\) is at most countable (and for curves of genus \(g >3\) finite and for a general curve empty). 
Now suppose that we have chosen \( \sL_1, \dots , \sL_k\) linearly independent pints in \( im(\phi)\).
Note that the span of the points \( \sL_1, \dots ,\sL_k\) is countable in \(\J(C)\) and hence countable in \( im(\phi)\). 
 Thus, the totality of the span of \( \sL_1, \dots ,\sL_k\) and the set of torsion points in \( im(\phi)\) is countable.  
 Furthermore, to avoid higher singularities if \( (x_1,z_1) \dots (x_k,z_k)\) are the pairs of points that are chosen then in the \(k+1\) pair of point \((x_{k+1},z_{k+1})\) such that \( x_{k+1} \neq x_i\) or \(z_i\) and similarly for \( z_{k+1}\). This amounts to removing finitely many curves in \(im(\phi)\).
 We can thus choose a point \( \sL_{k+1} \in im(\phi)\) linearly independent of the previous choices.
Identifying the corresponding pairs of points in \( C\) gives us a singular curve \(X\)  with ordinary double points and by Theorem \ref{highnode} \( \J(X)\) is anti-affine.
\end{proof}

\begin{remark}\label{5.4}
Let \( X\) be a curve of genus \( g\). If the number of cuspidal singularities in \( X\) is greater than \( g\) then by  \cite[Proposition 2.3]{Brion2} we have that \( \J(X)\) cannot be anti-affine group. This is in stark contrast to the nodal case  proved in the above theorem.
We note that in the cuspidal case the Jacobian is not semi-abelian, the linear part of  \( \J(X)\) is a unipotent group. 
\end{remark}

\section{Ordinary \(n\)-point Singularity} \label{Onp}
Let \(X\) be a curve with \( y\in X\) the only singular point of \( X\). We assume \( y\) to be an ordinary \(n\)-point singularity. 
Recall that an ordinary \( n\)-point singularity is obtained by identifying \(n\) distinct points on a smooth curve \( C\) with each local branch of \( C\) immersed in \(X\) so that they give rise to independent vectors in the Zariski tangent space . 
Let \( p: C\to X\) be the normalization and \( p^{-1}(y) =\{ x_1 , \cdots, x_n\}\) be the fibre over \( y\).
We will denote the Jacobian of \( X\) , i.e., the space of degree 0 line bundles on \( X\), as \( \J(X)\).
Pulling back a line bundle \( \Li \in \J(X)\)  by \(p\) gives a map \( \pi: \J(X)\to \J(C)\).
We will embed \( \J(X)\) in a \( \mathbb{P}^{n-1}\) bundle over \( \J(C)\).
This bundle  \( \tilde{\J}(X)\) is defined as follows:
For any line bundle \( \Li \in \J(C)\), consider \( n-1\) dimensional quotient \( Q\), of the direct sum of fibres \( \Li_{x_1} \oplus \cdots \oplus \Li_{x_n}\). 
The set of pairs \( (\Li , Q) \) together with the first projection is a \( \mathbb{P}^{n-1}\) bundle over \( \J(C)\).
Thus, \( \tilde{\J}(X)= \mathbb{P}(\mathcal{P}_{x_1} \oplus \cdots \oplus \mathcal{P}_{x_n})\), here \( \mathcal{P}\) denotes the Poincar\'e bundle on \( \J(C) \times C\) and \( \mathcal{P}_{x_j}= \mathcal{P}|_{\J(C) \times \{x_j\}},~ x_j \in p^{-1}(y)\), cf. \cite{Bhosle} for details.
The inclusion \(h: \J(X) \hookrightarrow \tilde{\J}(X)\) is as follows:
Given a line bundle \( \Li \in \J(X)\), we  have a natural exact sequence of sheaves:
\[
\begin{tikzcd}
0 \to \Li \to p_* p^* \Li \to Q \to 0.
\end{tikzcd} 
\]

Then the map \( h\) is given by \( h(\Li) = (p^*\Li , Q)\). Since the above sequence is natural, the map \( h\) is injective.
The quotient \( Q\) is supported at the ordinary \(n\)-point \( y\) and is naturally the quotient of \( p^*\Li(x_1) \oplus \cdots \oplus p^*\Li(x_n)\) by \( \Li_y\).

Following  the forthcoming  paper \cite{Bhosle1} we will now show that the image of the map \(h\) is the complement of the projections.
Let \( (A, \mathfrak{m} ) \) denote the local ring of the ordinary \(n \)-point singularity \( y \in X\) and let \( \bar{A}\) be the normalization. Let \( \mathfrak{m}_1, \cdots , \mathfrak{m}_n\) be the maximal ideals of \( \bar{A}\). 
Let \(k_i:= \bar{A}/\mathfrak{m}_i =k \). Then \( \bar{A}/ \mathfrak{m}\bar{A} \approx \bar{A}/\mathfrak{m}_1 \oplus \cdots \oplus \bar{A}/\mathfrak{m}_n \approx \mathop{\oplus}\limits_{i=1}^n k_i\)
and let \( \pi : \bar{A} \to k_1\oplus \cdots \oplus k_n\) denote the natural map.

\begin{lemma}
Let \( V\) be a line passing through the origin in \( k_1 \oplus \cdots \oplus k_n\). If all the \(n\) projections \( p_i: V \to k_i\) are isomorphisms then \( \pi^{-1} (V)\approx A\). 
\end{lemma}
\begin{proof} \renewcommand{\qedsymbol}{\(\blacksquare\)}
We will produce an automorphism \( f\) of \( \bar{A}\) such that the induced automorphism \( f_1\) of \( \bar{A}/\mathfrak{m}\bar{A}\) maps \( V\) onto the diagonal \( \Delta\) of \( k^n\).
The result follows as  \(\pi^{-1} (\Delta)=A\). \par
Let \( e\) be a generator of \( V\), and \( e = e_1 + \cdots e_n,~ e_i \in k_i\). Note since each \( p_i\) is an isomorphism  each \( e_i\) generates \( k_i\). Further, since \( \pi\) is surjective there is an element \( M \in \bar{A}\) such that \( \pi(M) = e\). Let \( f\) be the endomorphism of \( \bar{A}\) determined by \( M\). Let \( f_1: \bar{A}/\mathfrak{m}\bar{A}\to\bar{A}/\mathfrak{m}\bar{A}\) be the induced map. Then, note that \( f_1\) takes the canonical generator of \( k\) to \( e\). Thus \( M \notin \mathfrak{m}\) and hence is a unit. So, \( f\) is an isomorphism that lifts \( f_1\). It is clear that \( f_1(\Delta) =V\).

\end{proof}
The above lemma gives a criterion for the pushforward sheaf to be locally free.
We  have thus proved the following proposition.
\begin{propn}\label{semiabelian}\renewcommand{\qedsymbol}{\(blacksquare\)}
Let \( X\) be a projective curve with \( y\in X\) be the only singular point which we assume to be  an ordinary \(n\)-point singularity. Then \( \J(X)\) is a semiabelian variety and \( \J(X) = \mathbb{P}(\sP_{x_1} \oplus \cdots \oplus \sP_{x_n}) \setminus (s_{x_1} \bigcup \cdots  \bigcup s_{x_n}) \), where \( s_{x_j}\)s' are the divisors corresponding to the coordinate hyperplanes.
\end{propn}
\(\hfill \blacksquare\)

We are considering singular curves with ordinary singularities and we have shown that the Jacobian of such curves is semiabelian. 
We can ask the converse to this question, if the Jacobian of a singular curve  \( X\) is semiabelian then are the singularities of \( X\) ordinary?
We answer this question affirmatively. We first recall the definition of seminormality and see the equivalence of seminormality of \( X\),  \( X\) having ordinary singularities and \(\J(X)\) being semiabelian.

Following Swan \cite{Swan} we define seminormality as follows.
\begin{defn}[Seminormal ring] 
A commutative ring \( R\) is said to be seminormal  if whenever \( b, c \in R\) satisfy \( b^3 =c^2\) \ there exists \(a \in R\) such that \( a^2 =b\) and \( a^3 =c\).
\end{defn}
\begin{defn} [Seminormal scheme]
A scheme \( X\) is said to be seminormal if every \( x \in X\) has an affine open neighbourhood \( \Spec(R) = U \subset X\) such that \( R\) is seminormal.
\end{defn}
Following Koll\'ar \cite{Kollar}, we will give an alternative definition of seminormality.

\begin{defn} \label{KolSemi}

A noetherian scheme \( \) is said to be seminormal if for every morphism \( f: Y \to X\) that induces a homeomorphism (as a topological space ) and preserves the residue fields is an isomorphism.
\end{defn}

Thus, nodal curves are seminormal whereas cuspidal curves are not seminormal.
Geometrically, seminormality means that any non-normality (i.e., a gluing)  of a scheme will be as transverse as possible. We refer to \cite{Davis},\cite{Kollar},  \cite{Sal},  \cite{Swan}, \cite{Trav} for the properties and other variants of normality.

We will now prove the converse to Proposition \ref{semiabelian}. 

\begin{theorem}\label{seminorm}
The following are equivalent:
\begin{enumerate}
\item The curve \(X \) is seminormal.
\item The curve \(X \) has ordinary singularities. 
\item The Jacobian \(\J(X)\) is semiabelian.
\end{enumerate}
\end{theorem}
\begin{proof}\renewcommand{\qedsymbol}{\(\blacksquare\)}
The equivalences of 1 and 2 is   \cite[Corollary 4]{Davis}.
The proof that 3 implies 2 is  \cite[Chap. 9, Sect. 2 Corollary 12.a]{BLR}.

The proof that 2 implies 3 is the statement of Proposition  \ref{semiabelian}.
\end{proof}

We note that the result of the paper of Davis has been proved earlier by Bombieri \cite{Bomb}, however the paper is not yet digitized.

We will now give a criterion for a semiabelian Jacobian \( \J(X)\) to be anti-affine. 
We will first deal with the case of a single ordinary \(n\)-point singularity and then generalize to curves with multiple ordinary singularities.

Let \( X\) be a projective curve with  single singularity at \( y\) which we assume to be an ordinary \( n\)-point singularity.
Suppose \( p: C\to X\) be the normalization and \( p^{-1}(y)= \{x_1, \cdots , x_n\}\). 
To get a criterion for \( \J(X)\) to be anti-affine we follow the argument of Theorem \ref{node} and compute  \( \pi_*(\sO_{\J(X)})\).
First, we will twist the projective bundle by the line bundle \( \sP_{x_1}^{-1}\). For each \( i=2, \cdots , n\), denote \( \sP_{x_i}\otimes \sP_{x_1}^{-1}\) by \( \sL_i\). We have the following theorem.
 \begin{theorem}\label{npt}
 Let \( X, C\) and \( \{\sL_i\}_{i =1, \dots n}\) be as described in the previous paragraph. 
 Then the  Jacobian \( \J(X)\) is anti-affine if and only if the points \(\{ \sL_i\}_{i=2}^{n}\) are linearly independent in \( \J(C)\).
 \end{theorem}
 \begin{proof}\renewcommand{\qedsymbol}{\(\blacksquare\)}
 We will follow the proof of  Theorem \ref{node}.
 The regular functions on \( \J(X) \) are given by,  \( \sO( \J(X) )= \Ho^0( \J(C), ~\pi_*(\sO_{\J(X)}) \). We will give the description of \( \pi_*(\sO_{\J(X)}) \). 
Locally, by considering affine open sets of \( \J(C) \),  we have that \( \pi_*(\sO_{\J(X)}) \approx \mathop{\oplus}\limits_{\substack{k_i \in \mathbb{Z}\\\ i \in \{ 2,  \cdots n\}}} 
 \sL^{k_2}_2 \otimes \sL^{k_3}_3 \otimes \cdots \otimes  \sL^{k_{n}}_{n}\) and
  \( \Ho^0(\J(C), \pi_*(\sO_{ \J(X)} )) \approx \Ho^0 \big(\J(C), \mathop{\oplus}\limits_{\substack{k_i \in \mathbb{Z}\\\ i \in \{2,  \cdots n\}}} \sL^{k_2}_2 \otimes \sL^{k_3}_3 \otimes \cdots \otimes  \sL^{k_{n}}_{n}\big)\approx \mathop{\oplus}\limits_{\substack{k_i \in \mathbb{Z} \\\ i\in \{2,  \cdots n\}}}\Ho^0(\J(C),  \sL^{k_2}_2 \otimes \sL^{k_3}_3 \otimes \cdots \otimes  \sL^{k_{n}}_{n}  )\). 
 Therefore the group \( \J(X)\) is anti-affine if and only if the points \( \sL^{k_2}_2 \otimes \sL^{k_3}_3 \otimes \cdots \otimes  \sL^{k_{n}}_{n}\) are non trivial, 
i.e.,  the set of points \( \{\sL_2, \sL_3,\cdots ,\sL_{n}\}\) generate a free abelian group in \( \J(C)\).

 \end{proof}
 
\begin{remark}
The  torus extensions of an abelian variety \( A\) are in correspondence with group homomorphisms \(\Hom_{Grp} (\hat{T}, \hat{A})\) where \( \hat{T}\) is the character group of the torus \( T\) and \(\hat{A}\) is the dual abelian variety. Thus, for the principal \(T\)-bundle \( \J(X)\) to have the structure of an algebraic group twisting by a line bundle \( \sP_{x_1}^{-1}\) is essential , cf. Section 5.3 of \cite{Bri}. 
As the following lemma will show, the choice of the twist is immaterial.
 \end{remark}

\begin{lemma}\label{lemma}
Fix an index \( i_0 \in \{ 1, \cdots n\}\).
If the divisors \( \{(x_i-x_{i_0})\}_{\substack{ i\neq i_0\\i =1,\cdots, n}}\) are linearly independent in \( \J(C)\), then for any \( j_0 \neq i_0\)  the divisors \( \{(x_i-x_{j_0})\}_{\substack{ i\neq j_0\\i =1,\cdots, n}}\) are also linearly independent.
\end{lemma}
\begin{proof}\renewcommand{\qedsymbol}{\(\blacksquare\)}
For simplicity, we take \( i_0 =1\) and \( j_0 = 2\).

Suppose that \( \{(x_i-x_{1})\}_{\substack{ i\neq 1\\i =1,\cdots, n}}\) are linearly dependent. 
Then there are integers \( m_i, i =2,\dots n\) such that 
\begin{equation}\label{eqn1}
 \mathop{\sum}\limits_{i=2}^n m_i (x_i -x_1) =0 .
\end{equation}
We claim that \( \{(x_i-x_{2})\}_{\substack{ i\neq 2\\i =1,\cdots, n}}\) cannot be linearly independent. 
For this we add and subtract \( \mathop{\sum}\limits_{i=2}^n m_i(x_2)\) to the equation(\ref{eqn1}).
Thus, \(  \mathop{\sum}\limits_{i=2}^n m_i (x_i -x_1) -\mathop{\sum}\limits_{i=2}^n m_i(x_2) +\mathop{\sum}\limits_{i=2}^n m_i(x_2)\\
= -\mathop{\sum} \limits_{i=2}^n m_i(x_1 -x_2) + \mathop{\sum}\limits_{i=3}^n m_i(x_i -x_2) =0\). 
Thus the divisors \( \{(x_i-x_{2})\}_{\substack{ i\neq 2\\i =1,\cdots, n}}\) satisfy the linear equation 
\(\mathop{\sum}\limits_{i=1, i\neq 2}^n m_i' (x_i-x_2)\), with \(m_1' = \mathop{\sum}_{i= 2}^n m_i\) and for \(i\geq 3, m_i' = m_i\).
Thus, the divisors 
\( \{(x_i-x_{2})\}_{\substack{ i\neq 2\\i =1,\cdots, n}}\) are also not linearly independent.

\end{proof}

Let \( X\) be a curve with \( r\) ordinary-\( d_j\) singular points \( y_i, i=1, \cdots r\). 
Let \( \pi: C\to X\) be the normalization. 
For each \(i\), let \( \pi^{-1}(y_i) = \{ x_{i,1}, \cdots , x_{i,di}\}\) be the points in \( C\) corresponding to the singularity \( y_i\).
 As before, denote the Poincar\'e bundle on \( \J(C) \times C\) as \( \sP\).
Denote the vector  bundle  \( (\sP_{x_{i,2}} \otimes \sP_{x_{i,1}}^{-1})
\oplus \cdots \oplus 
 (\sP_{x_{i,d_i}} \oplus \sP_{x_{i,1}}^{-1}) \)  as \( \sV_i\) and the line bundle
\( (\sP_{x_{i,j}} \otimes \sP_{x_{i,1}}^{-1})\) as \( \sL_{ij}\) \( i= 1\dots , r\) and \(j = 2, \dots d_i\).
\begin{theorem} \label{6.5}
Let \( X\) be a projective curve as above.
Then the Jacobian \( \J(X) \hookrightarrow \bP(\sO\oplus \sV_{1})\times_{\J(C)} \cdots \times_{\J(C)} \bP(\sO\oplus \sV_{r})\). 
Furthermore \( \J(X)=\bP(\sO\oplus \sV_{1})\setminus \{\mathop{\bigcup}\limits_{i=1,\cdots d_1} s_i\}\times_{\J(C)} \cdots \times_{\J(C)} \bP(\sO\oplus \sV_{r})\setminus \{\mathop{\bigcup}\limits_{i=1,\cdots d_r} s_i\}\), where each \(s_i \) is a divisor  corresponding to a coordinate hyperplane.
Further, \( \J(X)\) is anti-affine if and only if the points \( \{\sL_{ij}\}_{\substack{i=1,\dots r \\\ j=2,\dots d_i}}\) are linearly independent in \( \J(C)\).
\end{theorem}
 \begin{proof}\renewcommand{\qedsymbol}{\( \blacksquare\)}
 The embedding of \( \J(X)\) in the projective bundle follows from Proposition \ref{semiabelian}.
 The proof of \( \J(X)\) is anti-affine is along the same line as Theorem \ref{highnode} and Theorem \ref{npt}.
 \end{proof}

Note that as before the choice of twisting is immaterial. Since the singularities are 'independent', i.e., the points in the normalization of each singularity are distinct from each other and  mixed tensor products do not appear in the computation of \( \pi_*(\sO_{\J(X)})\). Hence Lemma \ref{lemma} can be applied for the case of singular curves with multiple ordinary singularities.
 
 We will now give an analogue of Theorem \ref{choice} for curves with ordinary \(n\)-point singularity.
 
 \begin{propn}
 Let \( C\) be a smooth curve which is not hyperelliptic or bielliptic and let the genus \( g(C)\geq 3\). We fix the following data \( (n, d_1, \dots, d_n)\), \( n, d_i \in \mathbb{Z}^+\) and \( n\geq 1, d_i \geq 2\). Then there exists a curve \( X\) with exactly \( n\) singular points, each singular point is an ordinary \( d_i\)-point singularity and the Jacobian  \( \J(X)\) is anti-affine.
 \end{propn}
 \begin{proof}\renewcommand{\qedsymbol}{\( \blacksquare\)}
 We proceed by induction on \( n\) and \( d\).
 The base case is \((1,2)\), i.e., a nodal curve with a single double point. Note that this along with multiple ordinary double points is proved in Theorem \ref{choice}.
 
Before the general case we will show the case of a curve with single triple point, i.e., \( ( n=1, d=3)\)
Assume that we have made the choice of points \( x,w\) so that \( w-x\) is not torsion. 
We need to select a point \( z \in C\) so that \( z-x\) is not torsion and \(\{ (w-x), (z-x)\}\) are linearly independent. 
Note that such a \( z\) exists: for genus \( g(C)\geq 3\) the torsion divisors are finite (and for a general \( C\) empty cf. Theorem \ref{finale}) and the divisors that are linearly dependent on \( (w-x)\) are countable in \( \J(C)\) and hence countable in the surface \( im(\phi)\).
Identifying the points \(x,w,z\) we get a curve \( X\) with a triple point \(y\). 
By Theorem \ref{6.5} \( \J(X)\) is anti-affine.

For the general case: assume that we have chosen \(k-1\) singularities. 
Also assume that we have chosen points \( x_{k,1}, \dots , x_{k,l-1}\) so that \( \{x_{k,i}-x_{k,1}\}\) and the divisors from the other singular points are linearly independent and each \(x_{k,j}\) is distinct from the other choices of points. 
For genus \( g(C) \geq 3\) the torsion divisors of the form \( a-b\) are finite and the number of divisor that are linearly dependent on these divisors in \( \J(C)\) are countable and hence countable in the surface \( im(\phi)\).
Therefore we can find \( x_{k,l}\) so that \( x_{k,l}-x_{k,1}\) is linearly independent of the choices made.
This completes the induction. 
 \end{proof}


\begin{thebibliography}{9}

\bibitem{DH}
D. Abramovich and J. Harris, \textit{Abelian varieties and curves in \(W_d(C)\)}, Compositio Math. 78 (1991), no. 2, 227--238. 

\bibitem{Accola}
R. D. Accola,  \textit{Topics in the Theory of Riemann Surfaces},  Lecture Notes in Mathematics 1595,  Springer-Verlag.


\bibitem{Bhosle} 
U. Bhosle,  \textit{Generalised parabolic bundles and applications to torsion free sheaves on nodal curves},   
 Ark. Mat. 30 (1992),  no. 2,  187--215.
 
 \bibitem{Bhosle1}
 U. Bhosle, \textit{Representations of the fundamental group and Higgs bundles on singular integral curves},  to appear in Proc. Indian Acad. Sci. Math. Sci..

\bibitem{BP}
U. N. Bhosle,   A. J. Parameswaran \textit{On the Poincare formula and the Riemann singularity theorem over nodal curves},  Math. Ann. 342 (2008),  no. 4,  885--902. 

\bibitem{Bomb}
E. Bombieri, \textit{Seminormalita e singolarita ordinarie}, Sympos. math. 11, Algebra commut., Geometria, Convegni 1971/1972, 205-210. 1973. 

\bibitem{Borel}
A. Borel, \textit{Linear Algebraic Groups}, Second Enlarged Edition, Grad. Texts in Math. 126, Springer- Verlag, New York, 1991.

\bibitem{Bri}
M. Brion, P. Samuel, V. Uma, \textit{Lectures on the structure of algebraic groups and geometric applications}, CMI Lecture Series in Mathematics 1. New Delhi: Hindustan Book Agency; Chennai: Chennai Mathematical Institute (CMI), 2013.

\bibitem{Brion}
M. Brion, \textit{Some structure theorems for algebraic groups}, Proc. Sympos. Pure Math., 94, Amer. Math. Soc., Providence, RI, (2017), 53--125.

\bibitem{Brion1}
M. Brion, \textit{Commutative algebraic groups upto isogeny},  Doc. Math. 22 (2017), 679--725.

\bibitem{Brion2}
M. Brion,  \textit{Anti-affine algebraic groups},   J. Algebra 321 (2009),  no. 3,  934--952.


\bibitem{BLR}
S. Bosch, W. L\"utkebohmert, M. Raynaud, \textit{N\'eron Models}, Ergebnisse der Mathematik und ihrer Grenzgebiete, 3. Folge, 21. Berlin etc.: Springer-Verlag. x, 325 p. (1990)


\bibitem{Cyril}
C. D'Souza, \textit{Compactification of generalised Jacobians}, Proc. Indian Acad. Sci. Sect. A Math. Sci. 88 (1979), no. 5, 419--457.

\bibitem{Davis}
E. D. Davis, \textit{On the Geometric interpretation of Seminormality},
Proc. Am. Math. Soc. 68, 1-5 (1978).
\bibitem{HS}
J. Harris and J. Silverman,  \textit{Bielliptic Curves and Symmetric products},    Proc. Amer. Math. Soc. 112 (1991),  no. 2,  347--356. 

\bibitem{Hart}
R. Hartshorne,  \textit{Algebraic Geometry},  Graduate Texts in Mathematics,  Volume 52,  Springer-Verlag New York Heidelberg, 1977.


\bibitem{Kempf}
 G. R. Kempf, \textit{Torsion divisors on algebraic curves}, Pacific J. Math. 97 (1981), no. 2, 437--441.


\bibitem{Kollar}
J. Koll\'ar, \textit{Variants of Normality for Noetherian Schemes},  Pure Appl. Math. Q. 12 (2016), no. 1, 1--31.


\bibitem{Raynaud}
M. Raynaud,  \textit{Sous-vari\'{e}t\'{e}s d'une vari\'{e}t\'{e} ab\'{e}lienne et points de torsion},  Arithmetic and geometry,  Vol. \Rmnum{1},  Progr. Math.,  35,  Birkhäuser Boston,  Boston,  MA,  327--352 (1983).

\bibitem{Ray2}
 M. Raynaud, \textit{Courbes sur une vari\'et\'e abélienne et points de torsion}, Invent. Math. 71 (1983), no. 1, 207--233.

\bibitem{Rosen}
M. Rosenlicht, \textit{Toroidal algebraic groups},  Proc. Amer. Math. Soc. 12 (1961), 984--988. 

\bibitem{Sal}
P. Salmon \textit{Singolarità e gruppo di Picard}, 1969 Symposia Mathematica, Vol. II (INDAM, Rome, 1968) pp. 341–345.

\bibitem{Sancho}
C. Sancho de Salas, \textit{Grupos algebraicos y teoría de invariantes} Aportaciones Matemáticas: Textos, 16. Sociedad Matem\'atica Mexicana, M\'exico, 2001.


\bibitem{Sanchos}
C. Sancho de Salas and F. Sancho de Salas, \textit{ Principal bundles, quasi-abelian varieties and structure of algebraic groups},  J. Algebra 322 (2009), no. 8,  2751--2772.

\bibitem{Swan}
R. G. Swan, \textit{On seminormality}, J. Algebra 67 (1980), no. 1, 210--229. 

\bibitem{Trav}
C. Traverso, \textit{Seminormality and Picard group}, Ann. Sc. Norm. Super. Pisa, Sci. Fis. Mat., III. Ser. 24, (1970) 585--595.

\bibitem{Welt}
G. E. Welters, \textit{The surface C-C on Jacobi varieties and 2nd order theta functions}, Acta Math. 157 (1986), no. 1-2, 1--22.
\end{thebibliography}
\end{document}